\documentclass[12pt]{amsart}

\usepackage {amsfonts, amssymb, amscd, amsthm, amsmath, xypic}
\usepackage[dvips]{graphicx}
\usepackage{geometry}
\newtheorem{Thm}{Theorem}[section]
\newtheorem{theorem}[Thm]{Theorem}

\newcommand{\vep}{\varepsilon}
\renewcommand{\phi}{\varphi}
\newcommand{\co}{\mathbb{C}}
\newcommand{\N}{\mathbb{N}}

\newcommand{\D}{\mathbb{D}}
\newcommand{\T}{\mathbb{T}}
\newcommand{\tz}{\mathbb{T}}
\newcommand{\cp}{\mathbb{C_+}}
\newcommand{\R}{\mathbb{R}}
\newcommand{\rea}{{\rm Re}\,}

\title{On S. Grivaux' example of a hypercyclic rank one perturbation of a unitary operator.}
\author{Anton Baranov, Andrei Lishanskii}
\begin{document}
\sloppy
%\baselinestretched

%\large
%\normalsize

\address{ 
 Anton Baranov, 
\newline
Department of Mathematics and Mechanics,
St. Petersburg State University, 
\newline
St. Petersburg, Russia,
\newline
\phantom{x}\,\, and
\newline
National Research University  Higher School of Economics,
\newline
St. Petersburg, Russia,
\newline {\tt anton.d.baranov@gmail.com}
\newline\newline \phantom{x}\,\, Andrei Lishanskii,
\newline 
Department of Mathematics and Mechanics,
St. Petersburg State University, 
\newline 
St. Petersburg, Russia,
\newline
\phantom{x}\,\, and
\newline 
Chebyshev Laboratory,
St. Petersburg State University,
\newline 
St. Petersburg, Russia,
\newline {\tt lishanskiyaa@gmail.com}
\newline\newline \phantom{x}
}
\thanks{The authors were supported by the Chebyshev Laboratory 
(St. Petersburg State University) under RF Government grant 11.G34.31.0026 
and by JSC "Gazprom Neft". }

\begin{abstract}
Recently, Sophie Grivaux showed that there exists 
a rank one perturbation 
of a unitary operator in a Hilbert space which is hypercyclic.
We give a similar construction using a functional model
for rank one perturbations of singular unitary operators.
\end{abstract}

\maketitle

\section{Introduction}
A continuous linear operator $T$ in a Fr\'echet space $F$
is said to be {\it hypercyclic} if there exists a vector
$f\in F$ such that its orbit $\{T^n f\}_{n=0}^\infty$
is dense in $F$. In this case the vector $f$ is said
to be {\it hypercyclic for} $T$. First examples of 
hypercyclic operator go back to G.D. Birkhoff, G.R. McLane,
S. Rolewicz, while a systematic study of hypercyclicity phenomenon 
started in 1980-s with the thesis by C. Kitai and works by 
R.M. Gethner, G. Godefroy and J.H. Shapiro \cite{gesh, gosh}. We refer to the recent
monographs \cite{Bay-Mat, ge-per} for a detailed account of the theory.

Clearly, the identity operator is one of the most "nonhypercyclic"\, operators.
However, already in 1991 K.C. Chan and J.H. Shapiro \cite{chsh}
showed that there exists hypercyclic operators in a Hilbert space
of the form $I+K$, where the compact operator 
$K$ may belong to any Schatten class. It is clear that $I+R$ can not be hypercyclic when
$R$ is a finite  rank operator. Still, if we replace $I$ by a unitary operator, 
hypercyclicity  
is possible. In 2010 S. Shkarin \cite{shk} produced an example of a unitary operator
$U$ such that $U+R$ is hypercyclic for some rank two operator $R$. Shkarin asked
whether  $R$ can be taken to be of rank one. A positive answer was given by 
S. Grivaux \cite{gr1}:

\begin{theorem}[S. Grivaux, \cite{gr1}]
\label{main1}
There exists a unitary operator $U$ in the space $\ell^2$ 
and a rank one operator $R$ such that $U+R$ is hypercyclic.
\end{theorem}

The proof of this theorem is based on an ingenious elementary construction
involving a certain convergent inductive process as well as on the following
sufficient condition for hypercyclicity also obtained by Grivaux \cite{gr2}.
This result says that the operator is hypercyclic if there is a certain 
"continuous"\, family of eigenvectors with unimodular eigenvalues.

\begin{theorem}[S. Grivaux, \cite{gr2}]
\label{main2}
Let $X$ be a complex separable infinite-dimensional Banach space, and
let $T$ be a bounded operator on $X$. Suppose that there exists a sequence 
$\{u_n\}_{n\ge 1}$ of vectors in $X$ having the following properties:

{\rm (i)} $u_n$ is an eigenvector of $T$ associated to an eigenvalue $\lambda_n$ 
of $T$, with $|\lambda_n| = 1$ and $\lambda_n$ are all distinct;

{\rm (ii)} ${\rm span} \{u_n: n\ge 1\}$ is dense in $X$;

{\rm (iii)} for any $n \ge 1$ and any $\vep>0$, there exists 
$m\ne n$ such that $\|u_n -u_m\|<\vep$.
\\
Then $T$ is hypercyclic and even frequently hypercyclic.
\end{theorem}

The aim of the present paper is to give a proof of Theorem \ref{main1} 
by function theory methods. Our approach is based
on a functional model for rank one perturbations
of singular unitary operators. This model essentially goes back to 
a paper by V.V.~Kapustin \cite{kap}. In the form that we will use, this model 
appeared in \cite[Theorem 0.6]{bar-yak} in the context of rank one perturbations of
selfadjoint operators. This model translates
any rank one perturbation of a unitary operator to some concrete operator
in a space of analytic functions in the unit disk which is known  
as a star-invariant or model subspace (due to its role in yet another model -- 
that of B. Sz.-Nagy and C. Foias).

Let $H^2$ denote the standard Hardy space in the unit disk, and let 
$\theta$ be an inner function in the disk. The {\it model} (or {\it star-invariant}) 
{\it subspace} $K_\theta$ of $H^2$ is then defined as
$$
K_\theta=H^2\ominus \theta H^2.
$$
According to the famous Beurling theorem, any closed subspace of $H^2$ 
invariant with respect to the backward
shift in $H^2$ is of the form $K_\theta$. 
These subspaces play a distinguished role in operator theory (see,
e.g., \cite{Nikshift, Nik}) and in operator-related complex analysis.

The details  on the functional model for rank one perturbations
will be given in the next section. For the moment, let us mention only that in this 
model the family of the eigenvectors of a rank one perturbation have a very transparent
analytic meaning: they are either the families of reproducing kernels of 
$K_\theta$ or their biorthogonals. 

Now we state our main result which 
says that there exist model spaces with a certain continuous family of
vectors analogous to the properties of the vectors in Theorem \ref{main2}.
In view of the functional model (see Theorem \ref{rank-one-model}) 
and Theorem \ref{main2} this immediately
implies that there exists a unitary operator which has
a hypercyclic rank one perturbation.

\begin{theorem}
\label{main3}
There exists an inner function $\theta$ in the disk
such that $\theta(0) \ne 0$
and $\theta$ is analytically continuable across some nonempty open subarc of $\T$,
a function $\phi \in H^2 \setminus K_\theta$ and a sequence
$\lambda_n\in \T$ such that the functions
\begin{equation}
\label{main-eig}
f_n(z) = \frac{\phi(z)}{z-\lambda_n} \in K_\theta,
\end{equation}
the family $\{f_n\}$ is complete in $K_\theta$ and, 
for any $n \ge 1$ and any $\vep>0$, there exists 
$m\ne n$ such that $\|f_n -f_m\|<\vep$.
\end{theorem}

We do not by any means claim that our proof is essentially 
shorter than the original proof by Grivaux. However, 
we believe that the application of the model clarifies
the construction of the eigenvectors, since in the model space they 
have a special analytic structure: they are necessarily of 
the form \eqref{main-eig} for some $\phi\in H^2$.

As the original proof from \cite{gr1}, our construction is also inductive.
However, the parameters are chosen in a different way. In particular,
the eigenvalues of the operator $U+R$ will be chosen to be zeros of some Herglotz
function (interlacing with the spectrum of $U$). The properties
of Herglotz functions will play an important role in the construction.

\bigskip

%%%%%%%%%%%%%%%%%%%%%%%%%%%%%%%%%%%%%%%%%%%%%%%%%%%%%%
%%%%%%%%%%%%%%%%%%%%%%%%%%%%%%%%%%%%%%%%%%%%%%%%%%%%%%

\section{Preliminaries on the functional model for rank-one perturbations 
of unitary operators}

\subsection{Inner functions and Clark measures}
Recall that a function $\theta$ is said to be {\it inner} if it is analytic and bounded 
in the unit disk $\D$ and its nontangential boundary values satisfy
$|\theta| = 1$ a.e. with respect to the 
normalized Lebesgue measure $m$ on the unit circle $\mathbb{T}$.

Let $H^2 = H^2(\D)$ denote the {\it Hardy space} of the
unit disk $\mathbb{D}$, equipped with the standard norm
$\|\cdot\|_2 = \|\cdot\|_{L^2(m)}$. With each inner function 
$\theta$ we associate the model subspace $K_\theta
= H^2\ominus \theta H^2$. 
%In terms of the 
%nontangential boundary values $K_\theta$ may be defined as the space of all functions
%$f\in H^2$ such that $\bar{\zeta}\overline{f(\zeta)}\theta(\zeta) \in H^2$.

The {\it reproducing kernel} for $K_\theta$
corresponding to a point $\lambda\in\mathbb{D}$ is given by
$$
k_\lambda(z)=
\frac{1-\overline{\theta(\lambda)}\theta(z)}{1-\overline \lambda z}.
$$
Since functions in $K_\theta$ have more analyticity than general $H^2$
functions, there may exist reproducing kernels at boundary points.
In particular, one can consider $k_\lambda$, $\lambda \in I$, if $\theta$
(and, hence, any function in $K_\theta$) have an analytic continuation 
across the arc $I$. More generally, by the results of Ahern and Clark
%\cite{ac},
we have $k_\lambda \in K_\theta$ for $\lambda\in \tz$ if and only if
$|\theta'(\lambda)|<\infty$, the modulus of the angular derivative is finite.

Now we turn to Clark's construction
of orthogonal bases of reproducing kernels \cite{cl}.
For each $\alpha\in\mathbb{T}$, the function 
$\frac{\alpha+\theta}{\alpha-\theta}$ has positive real part in $\mathbb{D}$, and so
there exists a finite (singular) positive measure
$\mu^\alpha$ on $\tz$ such that
$$
\rea \frac{\alpha+\theta(z)}{\alpha-\theta(z)}=
\frac{1}{\pi} \int_\mathbb{T} \frac{1-|z|^2}{|\tau - z|^2}\,
d \mu^\alpha(\tau), \qquad z\in\mathbb{D}.
$$
Clark's theorem states that if, for some $\alpha$,
$\mu^\alpha$ is purely atomic,
i.e., if $\mu^\alpha = \sum_n \mu_n\, \delta_{\tau_n}$, $\tau_n \in \T$,
then $k_{\tau_n}\in {K_{\theta}}$ and 
the system $\{k_{\tau_n}\}$ is an orthogonal basis in ${K_{\theta}}$.
Note also that $\|k_{\tau_n}\|^2_2 = |\theta'(t_n)|= 2\mu_n^{-1}$.

Let $\mu = \mu^{1}$ be the Clark measure corresponding to $\alpha = 1$. 
For $c\in L^2(\mu)$, put 
\begin{equation}
\label{46}
(Vc)(z)=\big(1-\theta(z)\big)\int_\T
\frac{c(\tau) d\mu(\tau)}{1-\bar\tau z}, \qquad z\in \D.
\end{equation}
As Clark \cite{cl} has shown, $V$ is a unitary operator
from $L^2(\mu)$ onto $K_\theta$. Moreover, 
the nontangential boundary values of the function $Vc$ 
exist and coincide with $c$ $\mu$-a.e. \cite{polt}. 
In particular, if $\{k_{\tau_n}\}$ is an orthogonal basis
of reproducing kernels in $K_{\theta}$, then any $f\in K_\theta$
is of the form  $f(z)= \big(1-\theta(z)\big) \sum_n 
\frac{c_n\mu_n}{1-\bar\tau_n z}$, where $\sum_n |c_n|^2 \mu_n <\infty$.

%************************************************************

\subsection{The functional model}

Now we give the details of the functional model for rank one perturbations of unitary 
operators. 
Let $\theta$ be an inner function in the disk such that $\theta(0) \ne 0$.
For a function $f \in K_\theta$, the function $zf$ is not necessarily in $K_\theta$,
but it is easily seen that $zf = \gamma +h$, where $\gamma \in \co$ is a constant
and $h\in K_\theta$, and such decomposition is unique (moreover, $\gamma$ is a 
continuous functional of $f$). Now let $\phi \in H^2$
be a function such that
\begin{equation}
\label{main0}
\phi\notin K_\theta, \qquad \frac{\phi(z) -\phi(0)}{z} \in K_\theta.
\end{equation}
Then we may define the operator 
$T=T_{\theta,\phi}$ on $K_\theta$ by the formula 
\begin{equation}
\label{opt}
Tf := zf - \gamma_f \varphi, 
\end{equation}
where $\gamma_f$ is the unique complex number such that 
$zf - \gamma_f \varphi \in K_\theta$.
                                 
We are ready to present the functional model of rank one perturbations
which is analogous to \cite[Theorem 0.6]{bar-yak} (though the main ideas go back to 
\cite{kap}). Recall that a unitary operator is said to be {\it singular} if its 
spectral measure is singular with respect to the Lebesgue measure 
on $\T$. We also assume below that $U$ is cyclic, and so, up to a unitary 
equivalence, it is the operator of multiplication by $z$ in some space $L^2(\nu)$,
where $\nu$ is a finite Borel measure on $\T$. By $\sigma(U)$ we denote 
the spectrum of $U$. Finally, we denote by $\rho(\theta)$ the {\it boundary spectrum}
of an inner function $\theta$, that is, the complement of the union of all 
open arcs $I$ such that $\theta$ admits an analytic continuation       
across $I$.

\begin{theorem}[functional model]
\label{rank-one-model}
Let $U$ be a cyclic singular unitary operator such that $\sigma(U) \ne \T$.
Then for any rank one perturbation $U+R$ of $U$
there exist an inner function $\theta$ in the disk such that $\theta(0) \ne 0$
and $\rho(\theta) \ne \T$,
and a function $\phi\in H^2$ satisfying \eqref{main0} such that 
$U+R$ is unitary equivalent to the operator $T$ defined by \eqref{opt}. 
\smallskip

Conversely, any inner function $\theta$ 
such that $\theta(0) \ne 0$ and $\rho(\theta) \ne \T$,
and any function $\phi\in H^2$ satisfying \eqref{main0} correspond to 
some rank one perturbation $U+R$ 
of a cyclic singular unitary operator $U$ with $\sigma(U) \ne \T$.
\end{theorem}

It is obvious from the definition of the operator $T$ that if $\lambda\in \co$ is
an eigenvalue of $T$, then the corresponding eigenvector is given by 
$\frac{\phi(z)}{z-\lambda}$. Now, combining Theorems \ref{rank-one-model}
and \ref{main2} we see that the existence of a rank one perturbation 
follows from Theorem \ref{main3} (since the functions $f_n$ are eigenvectors 
of some rank one perturbations).

%************************************************************

\subsection{Model spaces in the upper half-plane}
\label{hp}
It will be more convenient to work with model spaces in the half-plane $\cp$
rather than in the disk. We can translate our problem to the half-plane setting
since the map $\tilde f(z) = \frac{1}{z+i} \cdot f\big(\frac{z-i}{z+i}\big)$
maps $H^2(\D)$ to $H^2(\cp)$ and the model space
$K_\theta$ in the disk to the model space 
$K_{\tilde \theta} = H^2(\cp) \ominus
\tilde\theta H^2(\cp)$ in $\cp$, 
where $\tilde \theta(z) = \theta\big(\frac{z-i}{z+i}\big)$.

The definition and the properties of the Clark measures for the model 
spaces in the upper half-plane are analogous to those in the disk.
Indeed, if $\theta$ is an inner function in $\cp$, then 
the function $\frac{\alpha+\theta}{\alpha-\theta}$ has positive real 
part in $\cp$, the measure $\mu$ on $\R$ from its Herglotz representation
is said to be a Clark measure for $K_\theta$, and, again, the embedding
of $K_\theta$ into $L^2(\mu)$ is a unitary operator (with one possible
exception where a linear term appears in the Herglotz representation; 
we will exclude this case in our construction). 

In particular, if the Clark measure corresponding to $\alpha =1$ is purely atomic, 
$\mu = \sum_n \mu_n \delta_{t_n}$, then any function $f\in K_\theta$ is of the form
\begin{equation}
\label{rep}
f(z) = \big(1-\theta(z)\big)\sum_n \frac{c_n\mu_n}{z- t_n}, \qquad
\sum_n |c_n|^2 \mu_n <\infty,
\end{equation}
and 
\begin{equation}
\label{clark}
\|f\|^2_{L^2(\R)} = 4\pi\sum_n |c_n|^2 \mu_n = \pi \sum_n |f(t_n)|^2\mu_n.
\end{equation}
We will often use this formula for the norm in what follows.

\bigskip

%%%%%%%%%%%%%%%%%%%%%%%%%%%%%%%%%%%%%%%%%%%%%%%%%%%%%%
%%%%%%%%%%%%%%%%%%%%%%%%%%%%%%%%%%%%%%%%%%%%%%%%%%%%%%

\section{Proof of Theorem \ref{main3}}
\label{proof}

\subsection{Plan of the solution}
By the above functional model, the problem is reduced to the construction
of functions $\theta$ and $\phi$ as in Theorem \ref{main3}. 
As explained in Subsection \ref{hp}, we can 
work on the real line and in the upper half-plane.
Thus, in what follows $\|\cdot\|_2$ stands for the usual $L^2$-norm on $\R$.

We will construct:
\begin{itemize}
\begin{item} 
a countable set $T = \{t_n\}_{n=1}^\infty$ on some interval
(say, in $[0,1]$);
\end{item}
\begin{item} 
a finite measure $\mu = \sum\limits_{n=1}^\infty \mu_n \delta_{t_n}$;
\end{item}
\begin{item} 
an inner function $\theta$ defined by
\begin{equation}
\label{t}
\frac{1 + \theta(z)}{1 - \theta(z)} := i \sum\limits_{n=1}^\infty \frac{\mu_n}{z - t_n};
\end{equation}
\end{item}
\begin{item} 
a function $\phi$ of the form
\begin{equation}
\label{p}
\varphi (z) := \big(1 - \theta (z)\big) \bigg[
\sum\limits_{n=1}^\infty \frac{c_n \mu_n}{t_n - z} + 1 \bigg], \qquad
\sum\limits_{n=1}^\infty |c_n|^2 \mu_n < \infty,
\end{equation} 
\end{item}
\end{itemize}
such that for some sequence  $\{\lambda_j\} \subset [0, 1],$ we have 
$\frac{\varphi}{z - \lambda_j} \in L^2(\mathbb{R})$ 
and, for any $j\ge 1$ and $\varepsilon > 0$ there exists $k\ne j$ such that
\begin{equation}
\label{st0}
\Big\| \frac{\varphi(z)}{z - \lambda_j} - \frac{\varphi(z)}{z - \lambda_k}
\Big\|_2 < \varepsilon.
\end{equation}
Note that in this construction $\mu$ is a Clark measure for $\theta$. 

As in Grivaux' paper \cite{gr1} we will proceed with the construction inductively. 
Namely, on the $N$-th step we construct
$t_1, \ldots, t_N \in [0, 1]$, $\mu_1, \ldots, \mu_N$ and functions 
$\theta_N$ defined by 
\begin{equation}
\label{tn}
i\frac{1 + \theta_N(z)}{1 - \theta_N(z)} 
:= \sum\limits_{n=1}^N \frac{\mu_n}{t_n - z},
\end{equation}
and  $\phi_N$,
\begin{equation}
\label{pn}
\varphi_N(z) := (1 - \theta_N(z))\bigg[1 + \sum\limits_{n=1}^N \frac{c_n \mu_n}{t_n - z}\bigg],
\end{equation}
where $c_n > 0$. The sequences $\mu_n$ and $c_n$ will be assumed to tend 
to zero very rapidly. 

It is a key idea of the construction that $c_n$ are taken to be {\it positive}. 
In this case the function $1 + \sum_{n=1}^N \frac{c_n \mu_n}{t_n - z}$ 
(which appears in the definition of $\phi_N$)
is a Herglotz function (has positive real part in $\cp$) and, therefore,
it has exactly $N$ zeros $\lambda_1^N, \lambda_2^N, \ldots \lambda_N^N$, 
interlacing with the points $t_1, \ldots, t_N$.
Choosing $c_n$ sufficiently small we can control the location of these points.

Note that $\theta_N$ is a finite Blaschke product and 
the model space $K_{\theta_N}$ is $N$-dimensional.
The measure $\sum_{n = 1}^N \mu_n \delta_{t_n}$ is a Clark measure for 
$K_{\theta_N}$.
Also, $\varphi_N \notin K_{\theta_N}$, since $1 - \theta_N \notin K_{\theta_N}$,
while 
$$
f_j^N(z) := \frac{\varphi_N(z)}{z - \lambda_j^N} \in K_{\theta_N}.
$$
Indeed, we have a representation of type \eqref{rep},
\begin{equation}
\label{rep1}
f_j^N(z) = \frac{\varphi_N(z) - \varphi_N(\lambda_j^N)}{z - \lambda_j^N} =
\big(1 - \theta_N(z)\big)  
\sum\limits_{n = 1}^N \frac{c_n \mu_n}{(\lambda_j^N - t_n)(z - t_n)}.
\end{equation}

Assume that  $t_n$, $\mu_n$ and $c_n$, $1\le n\le N_1$, are already chosen.
On $N$-th step we will add a point $t_{N}\in (0,1)$, its mass $\mu_{N}$ 
and a coefficient $c_{N}$ in the following order. First, we take
$t_{N}$ to be very close to some zero of the function $\phi_{N-1}$.
Then we choose $\mu_{N}$ to be so small that $\theta_{N}$ 
does not differ much from $\theta_{N-1}$ outside a small neighborhood
of $t_{N}$. Finally, we choose $c_{N}$ to be even much smaller than 
$\mu_{N}$ so that all zeros $\lambda_1^{N}, \ldots \lambda_N^{N}$ 
from generation $N$ almost coincide with the corresponding 
zeros $\lambda_1^{N-1}, \ldots \lambda_{N-1}^{N-1}$  from generation $N-1$, while
the zero $\lambda^{N}_{N}$ is very close to $t_{N}$.

Let us formally state what we need for the convergence.
\medskip

(I) We will choose $\mu_n>0$ and $c_n>0$ so that $\mu_n<2^{-n}$
and $c_n<2^{-n}$. These conditions already ensure the convergence 
of the functions \eqref{tn} and \eqref{pn} to 
\eqref{t} and \eqref{p} respectively. Also, we will require that
$|\lambda^{N-1}_j - \lambda^{N}_j| < 2^{-N}$, $j = 1, \ldots, N-1$
(in fact we will need much more, see \eqref{dist0} below).
\medskip

(II) Clearly, for each $N$, the functions $f_j^N$, $j=1, \ldots, N$, are linearly
independent, and so they form a basis in $K_{\theta_N}$. Let $A_N$ 
be a sort of a $\ell^1$ basis constant for $\{f_j^N\}$: for any 
$\{\alpha_j\}_{j=1}^N$, $\alpha_j\in \co$, 
$$
\sum_{j=1}^N|\alpha_j| \le A_N \Big\|\sum_{j=1}^N\alpha_j f_j^N\Big\|_2.
$$
Without loss of generality we assume that $A_N \ge 1$ and  
the sequence $A_N$ increases. Our second requirement then reads as follows:
\begin{equation}
\label{st}
\|f_j^{N} - f_j^{N-1}\|_2 < \frac{1}{2^{N+2} A_{N-1}}, \qquad j = 1, \ldots, N-1.
\end{equation}
\medskip

(III) Let $l(n)$ be a sequence of integers 
such that $l(n)<n$ and $l(n)$ takes every integer 
value infinitely many times (e.g., $1, 1, 2, 1, 2, 3, 1, 2, 3, 4, \ldots$).
To achieve property  \eqref{st0}, we will introduce the third requirement:
\begin{equation}
\label{st1}
\|f^{N}_{l(N)}- f^{N}_{N}\|_2 < 2^{-N-1}, \qquad N\in \N.
\end{equation}
\medskip

%**********************************************************

\subsection{Choice of the parameters}
\label{choice}
Assume that $t_n$, $\mu_n$ and $c_n$, 
$n=1, \ldots, N-1$, are already chosen.
First we choose the point $t_{N}$. Let 
$\varepsilon_{N}$ be some small positive number 
(namely, let $\vep_{N} \le 4^{-N-2}A_{N-1}^{-1}$, where $A_{N-1}$
is the constant from (II) above) and consider the equation
\begin{equation}
\label{eps}
1 + \sum\limits_{n = 1}^{N-1} \frac{c_n \mu_n}{t_n - x} = \varepsilon_{N}.
\end{equation}
Clearly, it has $N-1$ real zeros which 
depend continuously on $\vep_{N}$. Hence, 
they can be enumerated $x_1, \ldots x_{N-1}$ in such a way that 
$x_j \to \lambda_j^{N-1}$ when $\vep_{N} \to 0$. 

Let us take $\vep_{N}$ to be so small that  
if we take as $t_{N}$ the zero of \eqref{eps} which is the closest to the point 
$\lambda_{l(N)}^{N-1}$, then
\begin{equation}
\label{bbb}
|t_{N} - \lambda_{l(N)}^{N-1}| < 4^{-3N}\delta^3_{N-1},  
\end{equation}
where
$$
\delta_{N-1} = \inf_{1\le j, n\le N-1} |\lambda_j^{N-1} - t_n|.
$$

Now we put $\mu_{N} = \sqrt{\varepsilon_{N}}$ 
and define the inner function $\theta_{N}$ by \eqref{tn}.
One more additional restriction on the smallness 
of $\vep_{N}$ will be as follows:
\begin{equation}
\label{sm}
\Big\|\frac{\theta_{N}(z)-\theta_{N-1}(z)}{z-t_n}\Big\|_2 <\frac{1}{4^N A_{N-1}}, 
\qquad n=1, \ldots, N-1.
\end{equation}
This is also possible by continuous dependence of these norms from $\vep_{N}$.
More precisely, by the construction of $t_{N}$ we have
$t_{N} \to \lambda_{l(N)}^{N-1}$ as $\vep_{N} \to 0$, and so $t_{N}$ 
is separated from $t_n$, $1\le n\le N-1$.
Therefore, we can choose a small neighborhood $I$ of 
$\lambda_{l(N)}^{N-1}$ such that the integral of 
$\Big|\frac{\theta_{N}(t)-\theta_{N-1}(t)}{t-t_n}\Big|^2$ over 
$I$ is small. Once the interval $I$ is fixed, we can 
make the integral over $\R\setminus I$ to be arbitrarily small since
$\theta_{N}-\theta_{N-1}$ (as well as $\theta'_{N}-\theta'_{N-1}$)
tends to zero uniformly over $\R\setminus I$ as $\vep_{N} \to 0$.

Finally, let us choose $c_{N}$. 
Consider the equation $\phi_{N}(x) = 0$
which is equivalent to 
$$
1 + \sum\limits_{n=1}^{N-1} \frac{c_n \mu_n}{t_n - x} + 
\frac{c_{N} \mu_{N}}{t_{N} - x} =0.
$$
Again, a continuity argument shows that we may enumerate the $N$ zeros 
of this equation $\lambda_j^{N}$ (which should be 
thought of as functions of $c_{N}$) so that 
$\lambda_j^{N} \to \lambda_j^{N-1}$, $j=1, \ldots, N-1$, 
and $\lambda_{N}^{N} \to t_{N}$ as $c_{N} \to 0$.
Moreover, by the choice of $t_{N}$ as a solution of \eqref{eps}, we have
$$
1 + \sum\limits_{n = 1}^{N-1} \frac{c_n \mu_n}{t_n - \lambda_{N}^{N}} \to 
\varepsilon_{N}
$$
and $\lambda_{N}^{N} \to t_{N}$ as $c_{N} \to 0$.
Therefore, 
$$
\lambda_{N}^{N} - t_{N} \sim \frac{c_{N} \mu_{N}} {\vep_{N}}=
\frac{c_{N}} {\sqrt{\vep_{N}}}, \qquad c_{N} \to 0.
$$
Let us now choose $c_{N}$ to be so small that 
\begin{equation}
\label{dist0}
|\lambda_j^{N} -\lambda_j^{N-1}|< \frac{1}{4^{3N} A_{N-1}} \delta_{N-1}^3,
\qquad j=1, \ldots, N-1, 
\end{equation}
\begin{equation}
\label{dist}
\frac{c_{N}}{2\sqrt{\vep_{N}}} \le |\lambda_{N}^{N} - t_{N}| 
< |\lambda^{N}_{l(N)} - t_{N}|,
\end{equation}
and 
\begin{equation}
\label{dist2}
|\lambda_{N}^{N} - \lambda_{l(N)}^{N-1}| < 2^{-N}\delta^3_{N-1}.
\end{equation}
The latter estimate is possible by \eqref{bbb} and \eqref{dist0}.

%**********************************************************

\subsection{Proof of \eqref{st}}
To estimate the norm $\|f_j^{N} - f_j^{N-1}\|_2$, $1\le j\le N-1$, 
we write, using \eqref{rep1},
%$$
%\varphi_{N+1} (z) - \varphi_N (z) = 
%(1 - \theta_{N+1} (z)) \frac{c_{N+1} \mu_{N+1}}{t_{N+1} - z} 
%+ (\theta_N (z) - \theta_{N+1} (z)) \sum\limits_{n = 1}^N 
%\frac{c_n \mu_n}{t_n - z}.$$
%In the first summand of this expression we have by Lemma 1 that
%$$\frac{1 - \theta_{N+1}(z)}{t_{N+1} - z}|_{z = t_{N+1}} = \frac{2i}{\mu_{N+1}}.$$
%So $\|\frac{1 - \theta_{N+1}(z)}{t_{N+1} - z}\|^2_{L^2} = \|k_{t_{N+1}}\|^2_{L^2} = |k_{t_{N+1}}(t_{N+1})| = \frac{2}{\mu_{N+1}}$.
%In the second summand $\theta_N (z) - \theta_{N+1} (z)$ we can choose very small outside the neighborhood of the point $t_{N+1}$, because it is uniformly convergent to 0 out of the neighborhood of $t_{N+1}$, when $\mu_{N+1} \to 0$.
\begin{equation}
\label{sl}
\begin{aligned}
\frac{\varphi_{N-1}(z)}{z - \lambda_j^{N-1}} & - \frac{\varphi_{N}(z)}
{z - \lambda_j^{N}}  = \sum\limits_{n=1}^{N-1} 
\frac{c_n \mu_n (1 - \theta_{N-1}(z))}{(t_n - \lambda_j^{N-1})(z - t_n)} \\
& - \sum\limits_{n=1}^{N-1} \frac{c_n \mu_n(1 - \theta_{N})}
{(t_n - \lambda_j^{N})(z - t_n)} 
 - \frac{c_{N} \mu_{N} (1 - \theta_{N})}{(t_{N} - \lambda_j^{N})(z - t_{N})}.
\end{aligned}
\end{equation}
Denote the last summand in \eqref{sl} by $h$. Then $h\in K_{\theta_{N}}$
(it is a reproducing kernel up to a coefficient), and so 
$$
\|h\|_2^2 = \frac{|c_{N}|^2\mu_{N}}{|t_{N} - 
\lambda_{N}^{N}|^2} \le 4\mu_{N} \vep_{N}.
$$
The first two summands in \eqref{sl} may be rearranged into the sum of
$$
g_1(z) = \big(1 - \theta_{N}(z)\big)
\sum\limits_{n=1}^{N-1} \bigg(\frac{c_n \mu_n}{(t_n - \lambda_j^{N-1})(z - t_n)} - 
\frac{c_n \mu_n}{(t_n - \lambda_j^{N})(z - t_n)}\bigg)
$$
and 
$$
g_2(z) = \sum\limits_{n=1}^{N-1} c_n \mu_n \frac{\theta_{N}(z) - \theta_{N-1}(z)}
{(t_n - \lambda_j^{N-1})(z - t_n)}.
$$
We may compute the norm of $g_1\in K_{\theta_{N}}$
using Clark's theorem (see formula \eqref{clark}):
$$
\|g_1\|_2^2  = 
\sum\limits_{m=1}^{N-1} |g_1(t_m)|^2 \mu_m = 
\sum\limits_{m=1}^{N-1} \frac{|\lambda_j^{N} - \lambda_j^{N-1}|^2 |c_m|^2 \mu_m}
{|t_m - \lambda_j^{N-1}|^2|t_m - \lambda_j^{N}|^2}.
$$
Thus, $\|g_1\|_2$ does not exceed $4^{-N}A_{N-1}^{-1}$ by \eqref{dist0}.

Finally, $\|g_2\|_2 \le 4^{-N}A_{N-1}^{-1}$, and summing the above estimates 
we obtain \eqref{st}.

%**********************************************************

\subsection{Proof of \eqref{st1}}
To estimate the norm $\|f^{N}_{l(N)}- f^{N}_{N}\|_2$ we use again 
Clark's theorem:
\begin{equation}
\label{nnn}
\Big\|\frac{\varphi_{N}}{z - \lambda_{l(N)}^{N}} - 
\frac{\varphi_{N}}{z - \lambda_{N}^{N}}\Big\|^2_2 
= \sum\limits_{m = 1}^{N} \Big|\frac{\varphi_{N}(t_m)}{t_m - 
\lambda_{l(N)}^{N}} - \frac{\varphi_{N}(t_m)}{t_m - \lambda_{N}^{N}}\Big|^2 \mu_m.
\end{equation}
Let $L_{N} = |\lambda_{l(N)}^{N} - \lambda_{N}^{N}|$.
Then, by \eqref{dist0} and 
\eqref{dist2}, $|t_m - \lambda_j^{N}| > L_{N}^{1/3}$
for any $m = 1 \ldots N-1$ and $j=1, \ldots, N$. Now the first $N-1$ 
summands in \eqref{nnn} may be estimated as
$$
\sum\limits_{m=1}^{N-1} 
\frac{|\lambda_{l(N)}^{N} - 
\lambda_{N}^{N}|^2}{|t_m - \lambda_{l(N)}^{N}|^2 |t_m - \lambda_{N}^{N}|^2}
|\varphi_{N}(t_m)|^2 \mu_m   \le \sum\limits_{m=1}^{N-1} 
L_{N}^{2/3} |\varphi_{N}(t_m)|^2 \mu_m.
$$
Since by \eqref{dist0} and \eqref{dist2} $L_N \le 2^{-3N}\delta_{N-1}^3$ and
$|\phi_{N}(t_m)| = |\theta_{N}'(t_m)|c_m \mu_m = 2c_m$, we
conclude that the latter sum does not exceed $2^{-N-1}$. 

It remains to consider the summand with the number $N$. We have
$$
\bigg|\frac{\varphi_{N}(t_N)}{t_{N} - \lambda_{N}^{N}}\bigg| 
= \frac{2} {\mu_{N}} \cdot \frac{c_{N} \mu_{N}} {|t_{N} - \lambda_{N}^{N}|}
\le \frac{2\vep_{N}} {\mu_{N}} = 2\sqrt{\vep_{N}} \le 2^{-N-2}. 
$$
Here we used inequality \eqref{dist}. 
The estimate for the term $\Big|\frac{\varphi_{N}}{t_{N} - 
\lambda_{l(N)}^{N}}\Big|$ is analogous. 
Estimate \eqref{st1} is proved

%**********************************************************

\subsection{Convergence and completeness}

To ensure the pointwise  convergences $\theta_N(z) \to \theta(z)$
and $\phi_N(z)\to\phi(z)$, $z\in \cp$, it suffices to assume only that
$\sum_n \mu_n < \infty$ and $c_n \to 0$.

By the choice of the parameters in Subsection \ref{choice},  
the sequence  $\lambda_j^N$ converges to some point $\lambda_j$,
and it follows that $f_j^N$ converges to $\frac{\phi(z)}{z-\lambda_j}$ 
pointwise in $\cp$. On the other hand, by \eqref{st}, the sequence
$f_j^N \in H^2$ converges to some function $f_j$ in $H^2$ (recall that 
$H^2 = H^2(\cp)$ is the Hardy space in the upper half-plane). 
Since the convergence in $H^2$ implies the 
pointwise convergence in $\cp$, we conclude that 
$f_j(z) = \frac{\varphi(z)}{z - \lambda_j}$.

Next we prove that the family $\{f_j\}$ 
is complete in $K_\theta$. First of all, note that if $g\in K_\theta$,
$$ 
g(z) = \big(1 - \theta(z)\big) \sum\limits_{n=1}^\infty 
\frac{d_n\mu_n}{z - t_n}, \qquad 
g_N(z) = \big(1 - \theta_N(z)\big)\sum\limits_{n=1}^N \frac{d_n\mu_n}{z - t_n},
$$
then $\|g- g_N\|_2\to 0$, $N\to \infty$.
Indeed, 
$$
g(z)-g_N(z) = \big(\theta_N(z)-\theta(z)\big)\sum\limits_{n=1}^N 
\frac{d_n\mu_n}{z - t_n} + 
\big(1 - \theta(z)\big) \sum\limits_{n = N+1}^\infty \frac{d_n}{z - t_n}.
$$
The norm of the second sum by Clark's theorem (see \eqref{clark}) equals 
$\sum_{m=N+1}^\infty |d_m|^2\mu_m$, which obviously goes to zero, $N\to \infty$,
while the norm of the first sum is small by the assumption \eqref{sm}.

Thus, we constructed a sequence $g_N\in K_{\theta_N} $
such that $g_N \to g$ in $L^2(\R)$.
It remains to approximate functions $g_N$
by linear combinations of $f_j$.
The method is borrowed from \cite{gr1}. 
Since $\{f_j^N\}$ is a basis in $K_{\theta_N}$, we may write 
$g_N = \sum_{j = 1}^N \alpha_j f_j^N$.
Then, making use of \eqref{st}, we get 
$$
\begin{aligned}
\|g - \sum\limits_{j = 1}^N \alpha_j f_j\|_2 & \leq 
\sum\limits_{j = 1}^N |\alpha_j| 
\sum_{k=N}^{\infty} \|f_j^k - f_j^{k+1}\|_2  \\
& \le \frac{1}{2^N A_N}\sum_{j=1}^N|\alpha_j|\le 2^{-N}\|g_N\|_2,
\end{aligned}
$$
which goes to $0$ as $N\to \infty$. Completeness of the
family $\{f_j\}$ is proved.

%**********************************************************

\subsection{End of the proof of Theorem \ref{main3}.}
Let us complete the proof of Theorem \ref{main3}. We have constructed 
a complete sequence $\{f_j\} = \{ \frac{\phi}{z-\lambda_j}\}$
in $K_\theta$. It remains to verify that $\{f_j\}$ has the property
\eqref{st0}.  Let $j\ge 1$ and $\vep$ be given. Choose $N$ such that
$l(N) = j$ and $2^{-N} <\vep$, which is possible by the definition of 
the sequence $l(n)$. Then, by \eqref{st1}, $\|f^N_j - f^N_N\|_2 <2^{-N-1}$.
Also,  by \eqref{st1},
$$
\|f_N^N - f_N\|_2 \le \sum_{k= N}^\infty \|f_N^k - f_N^{k+1}\|_2 \le 2^{-N-2}
$$
and, analogously, 
$\|f_j^N - f_j\|_2 \le 2^{-N-1}$. Combining these estimates we 
obtain $\|f_j - f_N\|_2 \le 2^{-N}$. 
\qed

\bigskip

%%%%%%%%%%%%%%%%%%%%%%%%%%%%%%%%%%%%%%%%%%%%%%%%%%%%%%
%%%%%%%%%%%%%%%%%%%%%%%%%%%%%%%%%%%%%%%%%%%%%%%%%%%%%%

\section{Concluding remarks}

The unitary operator $U$ in our construction
(as well as in \cite{gr1}) is of a  very special form.
Note that if $t_n\in\R$ is the sequence constructed in
Section \ref{proof}, then the spectrum of $U$ is given by $\{\tau_n\}$, 
$\tau_n = \frac{t_n-i}{t_n+i}$. The points $t_n$ are chosen inductively so that
$t_{N}$ be close to some of the points $t_n$, $1 \le n \le N-1$, and so 
the set $\{t_n\}$ has a certain self-similarity. It seems to be a natural
question, which unitary operators have hypercyclic 
rank one perturbations.
\medskip
\\
{\bf Problem 1.} To describe cyclic unitary operators $U$
such that $U+R$ is hypercyclic for some rank one operator $R$.
\medskip

In particular, it is not clear whether the spectrum $\sigma(U)$
can have nonempty interiour or positive measure.
\medskip
\\
{\bf Problem 2.} Does there exist a 
cyclic unitary operator $U$
such that $U+R$ is hypercyclic for some rank one operator $R$ and $\sigma(U) = \T$.
\medskip

As in \cite{gr1}, in the present paper the spectral measure of 
$U$ is purely atomic.
\medskip
\\
{\bf Problem 3.} Construct a unitary operator $U$,
whose spectral measure is an absolutely continuous 
or a continous (i.e., without point masses) singular  measure on $\T$,
such that $U+R$ is hypercyclic for some rank one operator $R$.
\medskip

Note that the functional model applies to rank one perturbations of an arbitrary
cyclic unitary operator whose spectral measure is singular. Therefore, one can hope
to obtain further information about hypercyclic rank one perturbations of unitary 
operators using this model.

\end{document}